\documentclass[a4paper, 12pt]{amsart}
\usepackage{amsmath, amssymb, amsthm, graphicx, stmaryrd}
\usepackage{mathrsfs}
\usepackage{comment}
\usepackage[all]{xy}
\usepackage{tikz}
\usetikzlibrary{decorations.markings}

\usepackage[top=40mm,bottom=40mm,left=30mm,right=30mm]{geometry}

\theoremstyle{definition}
\newtheorem{theorem}{Theorem}[section]
\newtheorem{definition}[theorem]{\rm Definition}
\newtheorem{lemma}[theorem]{Lemma}
\newtheorem{proposition}[theorem]{Proposition}
\newtheorem{corollary}[theorem]{Corollary}
\newtheorem{remark}[theorem]{\rm Remark}
\newtheorem{example}[theorem]{\rm Example}

\newtheorem*{Atheorem}{Theorem A}
\newtheorem*{Btheorem}{Theorem B}

\newtheorem*{pf}{\rm Proof}

\makeatletter
\@addtoreset{equation}{section} 
\makeatother

\DeclareMathOperator{\Diff}{Diff}

\DeclareMathOperator{\id}{id}

\DeclareMathOperator{\Symp}{Symp}
\DeclareMathOperator{\rel}{rel}
\DeclareMathOperator{\flux}{Flux}
\DeclareMathOperator{\Ker}{Ker}
\DeclareMathOperator{\Cal}{Cal}

\title[The flux homomorphism and central extensions]{The flux homomorphism and central extensions\\
of diffeomorphism groups}
\author{Shuhei Maruyama}
\address{Graduate School of Mathematics,
Nagoya University, Japan}
\email{m17037h@math.nagoya-u.ac.jp}

\begin{document}

\begin{abstract}
  Let $D$ be a closed unit disk and
  $G_{\rel}$ the group of symplectomorphisms preserving
  the origin and the boundary $\partial D$ pointwise.
  We consider the flux homomorphism on $G_{\rel}$
  and construct a central $\mathbb{R}$-extension called
  the flux extension.
  We determine the Euler class of this extension
  and investigate the relation among the extension,
  the group $2$-cocycle defined by
  Ismagilov, Losik, and Michor 
  and the Calabi invariant of $D$.
\end{abstract}

\maketitle

\section{Introduction}

\footnote[0]{2010 Mathematics Subject Classification. 55R40, 37E30}

Let $D = \{ (x, y) \in \mathbb{R}^2 \mid x^2 + y^2 \leq 1 \}$
be a closed unit disk
with a symplectic form
$\omega = dx \wedge dy$.
Denote by
$G = \{ g \in \Diff(D)
\mid
\omega^g = \omega, g(0,0) = (0,0) \}$
the group consisting of
symplectomorphisms on $D$ that preserve the origin
and set
$G_{\rel}
=
\{ g \in G \mid g|_{\partial D} = \id_{\partial D} \}$.
Then there is the following exact sequence:
\[
  1 \longrightarrow G_{\rel} \longrightarrow
  G \longrightarrow \Diff_{+}(S^1)
  \longrightarrow 1,
\]
where $\Diff_{+}(S^1)$ is the group of orientation preserving
diffeomorphisms on the unit circle $S^1 = \partial D$.
On the group $G_{\rel}$,
there is an $\mathbb{R}$-valued homomorphism
\[
  \flux_{\mathbb{R}} : G_{\rel} \to \mathbb{R}
  \ ; g \mapsto \int_{\gamma} g^*\eta - \eta
\]
called a flux homomorphism,
which is a composing of the integration and the
ordinary flux homomorphism
(see \cite{mcduff_salamon} for the ordinary flux
homomorphism).
Dividing by the kernel $K = \Ker \flux_{\mathbb{R}}$, we obtain
the following central $\mathbb{R}$-extension:
\[
  0 \longrightarrow \mathbb{R} \longrightarrow
  G/K \longrightarrow \Diff_{+}(S^1)
  \longrightarrow 1,
\]
which we call the flux extension.
Denote by $e(G/K)$
the Euler class of the flux extension,
namely the cohomology class in $H^2(\Diff_+(S^1);\mathbb{R})$
corresponding to the given extension.
There is another cohomology class $e_{\mathbb{R}}$
in $H^2(\Diff_{+}(S^1);\mathbb{R})$,
which is a class corresponding to a universal covering
space of $\Diff_{+}(S^1)$.
The following theorem clarifies the relation
between $e(G/K)$ and $e_{\mathbb{R}}$:
\begin{Atheorem}[Theorem \ref{main_thm_1}]
  The class $e(G/K)$ is equal to $\pi e_{\mathbb{R}}$.
\end{Atheorem}


For a symplectic manifold $M$ endowed with an
exact symplectic form $\omega = d\eta$,
there is a $2$-cocycle $C_{\eta, x_0}$ on $\Symp(M)$
defined by Ismagilov, Losik, and Michor\,\cite{IsLoMi06},
which we call the ILM cocycle.
Let us consider the ILM cocycle on disk.
It turns out that the ILM cocycle is related to
the flux homomorphism.
Let $\tau: G \to \mathbb{R}$ be the map defined
by the same as flux homomorphism,
that is, $\tau(g) = \int_{\gamma}g^*\eta - \eta$.
Then the ILM cocycle
is equal to the coboundary $-\delta \tau$
on $G$.

The ILM cocycle on disk also relates to the
Calabi invariant.
Set $H = \Symp(D)$ and $H_{\rel} = \{ g \in H \mid
g|_{\partial D} = \id \}$.
The Calabi invariant $\Cal : H_{\rel} \to \mathbb{R}$
is defined as
\[
  \Cal (g) = \int_{\partial D} \eta \cup \delta \eta(g)
  = \int_{\partial D} g^*\eta \wedge \eta.
\]
Let $\tau_0: H \to \mathbb{R}$ be
the function defined by the same formula of $\Cal$.
Then it turns out that the ILM cocycle coincides with
$\delta \tau_0$.
From this, we have;

\begin{Btheorem}[Theorem \ref{main_thm_3}]
  The ILM cocycle
  $C_{\eta,x_0}$
  is basic, that is,
  there exists a cocycle $\chi$ on $\Diff_{+}(S^1)$
  such that $C_{\eta, x_0} = p^*\chi$ with
  $p:H \to \Diff_{+}(S^1)$ the restriction.
  Furthermore, the cohomology class $[\chi]$ is
  equal to
  $\pi e_{\mathbb{R}}$.
\end{Btheorem}
\noindent

The present paper is organized as follows.
In section $2$, we briefly recall the Euler
class of a central extension.
We descrive a cocycle representing the Euler
class in terms of a connection cochain and its curvature.
In section $3$, we construct the flux extension
and we prove Theorem A.
In section $4$, we introduce the ILM cocycle.
Finally, in section $5$,
we discuss the relation between
the Calabi extension and the ILM cocycle,
and prove Theorem B.


\section{Central extensions and the Euler class}
Let $\Gamma$ be a group and $A$ a right $\Gamma$-module.
We denote the action by $a^g$ for $a \in A, g \in \Gamma$.
For a non-negative integer $p$,
a $p$-cochain is an arbitrary map
$c:\Gamma^p \to A$
where $\Gamma^p$ denote a $p$-fold product group.
Set
$C^p(\Gamma;A) = \{ c : \Gamma^p \to A \mid \text{$p$-cochains} \}$.
These cochains gives rise to a cochain complex
$(C^{\bullet}(\Gamma, A),\delta)$
and its cohomology $H^*(\Gamma;A)$ is called
a group cohomology
(see \cite{brown82} for more details).

If $A$ is a trivial $\Gamma$-module,
the group cohomology $H^2(\Gamma;A)$
has a description in terms of
central extensions of groups.
Recall that a \textit{central $A$-extension of $\Gamma$} is
a exact sequence of groups
\[
  0 \longrightarrow A \overset{i}{\longrightarrow} E
  \overset{p}{\longrightarrow} \Gamma \longrightarrow 1
\]
such that the image $i(A)$ is contained
in the center of $E$.
It is known that the each equivalence class of central $A$-extension
corresponds to a second cohomology class in $H^2(\Gamma:A)$
(see \cite{brown82}).
For a central $A$-extension $E$,
the corresponding cohomology class
$e(E) \in H^2(\Gamma;A)$
is called the
\textit{Euler class of the central $A$-extension $E$}.

To investigate the Euler class at the cochain level,
we introduce the notions called a connection cochain
and its curvature.
This is similar to the Chern-Weil theory
defining the Euler class of principal bundles.

\begin{definition}
  Let $E$ be a central $A$-extension of $\Gamma$.
  A $1$-cochain $\tau \in C^1(E;A)$ is called a
  \textit{connection cochain of $E$}
  if $\tau$ satisfies
  \[
    \tau(ga) = \tau(g) + a
  \]
  for any $g \in \Gamma, a \in A$.
  The coboundary $\delta \tau \in C^2(E;A)$ is
  called a \textit{curvature of} $\tau$.
\end{definition}

For a connection $\tau$,
its coboundary $\delta \tau$ is called a curvature.
A curvature is basic, that is, there exists a unique $2$-cocycle
$\sigma \in C^2(\Gamma;A)$
such that the pullback of $\sigma$ to $E$ coincides with
$\delta \tau$.
We call this cocycle $\sigma$ a basic cocycle of $\delta \tau$.
Then the cohomology class $[-\sigma]$ is equal to
the Euler class $e(E)$
(see \cite{Mori16} for more details).

Let $B$ be an abelian group and
$\iota:A \to B$ a homomorphism.
A $1$-cochain $\tau_B \in C^1(E;B)$
is called
a \textit{($B$-valued) connection cochain} if
$\tau_B$
satisfies $\tau_B(gb) = \tau_B(g) + b$
for every $g \in \Gamma, b \in B$.
The curvature $\delta \tau_B$ is basic and
we denote by $\sigma_B \in C^2(\Gamma;B)$ the basic cocycle.
Then the class $[-\sigma_B]$
corresponds to the image of Euler class $e(E)$
with respect to a natural homomorphism
$H^2(\Gamma;A)\to H^2(\Gamma;B)$.

\begin{example}
  Let us consider the Euler class
  $e_{\mathbb{R}} \in H^2(\Diff_{+}(S^1);\mathbb{R})$
  of a central extension
  $0 \to \mathbb{Z} \to
  \widetilde{\Diff_{+}(S^1)} \to \Diff_{+}(S^1)
  \to 1$.
  A $1$-cochain $\tau$ on $\widetilde{\Diff_{+}(S^1)}$
  is defined by
  $\tau (\varphi) = \varphi(0)/2\pi \in \mathbb{R}$
  where an element $\varphi$ in $\widetilde{\Diff_{+}(S^1)}$
  is considered as a diffeomorphism on $\mathbb{R}$
  satisfying $\varphi(x + 2\pi) = \varphi(x) + 2\pi$
  for any $x \in \mathbb{R}$.
  This cochain $\tau$ is a ($\mathbb{R}$-valued) connection cochain
  and the basic cocycle $\chi$ is
  \begin{equation}\label{univ_Euler_cocycle}
    \chi (\mu ,\nu)
    =
    -\delta \tau(\varphi, \psi)
    =
    \frac{\varphi\psi(0) - \varphi(0) - \psi(0)}{2\pi}
  \end{equation}
  where $\varphi,\psi \in \widetilde{\Diff_{+}(S^1)}$
  are lifts of $\mu ,\nu \in \Diff_{+}(S^1)$ respectively.
  This cochain $\chi$ is one of a cocycle of Euler class
  $e_{\mathbb{R}}$.
\end{example}

\section{The flux extension}
Next we introduce another central $\mathbb{R}$-extension
of $\Diff_{+}(S^1)$.
Let $G$ be the group of symplectomorphisms
preserving the origin $o = (0,0) \in \mathbb{R}^2$;
\[
  G = \{ g \in \Diff(D)
  \mid
  \omega^g = \omega, g(o) = o \}.
\]
We denote by $G_{\rel}$ the
subgroup of $G$ consisting of
symplectomorphisms which are the identity
on boundary.
Then there is
the following exact sequence;
\[
  1 \longrightarrow G_{\rel} \longrightarrow
  G \longrightarrow \Diff_{+}(S^1)
  \longrightarrow 1.
\]

We set $\eta = (xdy - ydx)/2$,
which satisfies $d\eta = \omega$.
Let $\gamma:[0,1] \to D$ be the singular $1$-chain defined by
\[
  \gamma(t) = (t,0) \in D
\]
for $t \in [0,1]$.

\begin{definition}
  The \textit{flux homomorphism}
  $\flux_{\mathbb{R}} \in C^1(G_{\rel};\mathbb{R})$
  is defined as
  \begin{equation}\label{def_eq_flux}
    \flux_{\mathbb{R}} (g)
    =
    -\int_{\gamma} \delta \eta (g)
    =
    \int_{\gamma} \eta^g - \eta
  \end{equation}
  for $g \in G_{\rel}$.
\end{definition}

\begin{proposition}
  The flux homomorphism is surjective.
\end{proposition}

\begin{pf}
  Let $\alpha : [0,1] \to \mathbb{R}$
  be the non-negative $C^{\infty}$ function
  sach that
  $\alpha(x) = 0$ on a neighborhood of $\{ 0,1 \}$ and
  $\alpha(1/2) = 1$.
  For $s>0$, define the symplectomorphism $h_s$
  by
  $h_s (r, \theta) = (r, \theta + s\alpha(r))$
  where $(r,\theta) \in D$ is the polar coordinates.
  Then we have $\flux_{\mathbb{R}}(h_1) >0$ and
  $\flux_{\mathbb{R}}(h_s) = s\flux_{\mathbb{R}}(h_1)$
  and the surjectivity follows.
\end{pf}

The right-hand side of formula (\ref{def_eq_flux})
also defines a $1$-cochain
$\tau: G \to \mathbb{R}$;
\begin{equation}\label{def_eq_tau}
  \tau(g) = -\int_{\gamma} \delta \eta (g).
\end{equation}

\begin{proposition}\label{prop_delta_tau}
  For $g,h \in G$,
  the following holds;
  \begin{equation}\label{eq_delta_tau}
    -\delta \tau (g,h)
    =
    \frac{1}{2}(\varphi \psi(0) - \varphi(0) - \psi(0))
  \end{equation}
  where
  $\varphi,\psi \in \widetilde{\Diff_{+}(S^1)}$ are lifts of
  the restrictions
  $g|_{\partial D},h|_{\partial D} \in \Diff_{+}(S^1)$
  respectively.
\end{proposition}

\begin{pf}
  The coboundary of $\tau$ is
  \begin{equation}\label{delta_tau_int}
    -\delta \tau (g,h)
    =
    \int_{\gamma - h\gamma}
    \delta \eta (g).
  \end{equation}
  By the Stokes formula and the exactness of
  $\delta \eta (g)$,
  the integration (\ref{delta_tau_int}) depends only on
  the endpoints of $\gamma - h\gamma$.
  So we have
  \begin{equation}\label{eq_delta_tau_and_ILM}
    -\delta \tau (g,h)
    =
    \int_{h(x_0)}^{x_0} \delta \eta (g)
    =
    - \int_{x_0}^{h(x_0)} \delta \eta (g)
  \end{equation}
  where $x_0 = (1, 0)$.
  Using the polar coordinates $(r,\theta)$,
  we put
  \[
    g(r,\theta) = (g_1(r,\theta),g_2(r,\theta)) \in D.
  \]
  Then the restriction of the integrand $\delta\eta(g)$
  to the boundary $\partial D$ becomes
  $\frac{1}{2}(d\theta - dg_2)$.
  Since $\varphi$ is a lift of
  $g|_{\partial D} (\theta) = g_2(1,\theta) \in \Diff_{+}(S^1)$,
  we obtain
  \begin{align*}
    -\delta \tau(g,h)
    &=
    - \frac{1}{2} \int_{x_0}^{h(x_0)}d\theta - dg_2
    =
    - \frac{1}{2} \int_{0}^{\psi(0)}
    dx - \frac{\partial \varphi}{\partial x}(x)dx\\
    &=
    \frac{1}{2}(\varphi\psi(0) - \varphi(0) - \psi(0)).
  \end{align*}
\end{pf}

\begin{remark}
  The equality (\ref{eq_delta_tau})
  implies that the cochain
  $\tau$
  is a quasi-morphism, that is,
  $\delta \tau$ is a bounded function.
  In fact,
  the absolute value
  $\frac{1}{2}|\varphi \psi (0) - \varphi (0) - \psi(0)|$
  is bounded by $\pi$
  (see \cite{ghys01} for more details).
\end{remark}

Next we show the formulas which are similar to
Kotschick-Morita\cite[Lemma $6$]{KoMo}.

\begin{corollary}\label{cor_of_flux_and_tau}
  For $g \in G, h \in G_{\rel}$,
  the following hold:
  \begin{enumerate}
    \item $\tau(gh) = \tau(g) + \flux_{\mathbb{R}}(h)$ and
    $\tau(hg) = \tau(g) + \flux_{\mathbb{R}}(h)$.
    \item The map $\flux_{\mathbb{R}} : G_{\rel} \to \mathbb{R}$
    is a homomorphism satisfying
    \[
      \flux_{\mathbb{R}} (ghg^{-1}) = \flux_{\mathbb{R}} (h).
    \]
    \item The kernel of the flux homomorphism
    $K = \Ker \flux_{\mathbb{R}}$ is a normal subgroup of $G$.
  \end{enumerate}
\end{corollary}

\begin{pf}
    As in the Proposition \ref{prop_delta_tau},
    we take lifts $\varphi$ and $\psi$ in
    $\widetilde{\Diff_{+}(S^1)}$
    of $g$ and $h$.
    Since $h \in G_{\rel}$, we have
    $\psi = T^n$ for some integer $n \in \mathbb{Z}$
    where $T:\mathbb{R} \to \mathbb{R}$ is a translation $T(x) = x + 2\pi$.
    Thus, we obtain
    \[
      \delta \tau (g,h)
      =
      \frac{1}{2}(\varphi(0)- \varphi T^n(0) + T^n(0))
      =
      \frac{1}{2}(\varphi(0) - \varphi(0) - 2n\pi + 2n\pi)
      =
      0
    \]
    and
    \[
      \delta \tau (h,g)
      =
      \frac{1}{2}(T^n(0) - T^n \varphi(0) + \varphi(0))
      =
      \frac{1}{2}(2n\pi - \varphi(0) - 2n\pi + \varphi(0))
      =
      0.
    \]
    On the other hand, we have
    \[
      \delta \tau(g,h)
      =
      \flux_{\mathbb{R}} (h) - \tau(gh) + \tau(g)
      \hspace{2mm}
      \text{ and }
      \hspace{2mm}
      \delta \tau(h,g)
      =
      \tau(g) - \tau(hg) + \flux_{\mathbb{R}}(h),
    \]
    which proves i).
    The equalities in i) prove ii) and the equality in
    ii) proves iii).
\end{pf}

According to the Corollary \ref{cor_of_flux_and_tau},
we have the exact sequence
\[
  0 \longrightarrow \mathbb{R} \longrightarrow
  G/K \longrightarrow \Diff_{+}(S^1)
  \longrightarrow 1
\]
where $\mathbb{R}$ is identified with
$G_{\rel}/K$.

\begin{proposition}
  The sequence
  \begin{equation}\label{rel_flux_ext}
    0 \longrightarrow \mathbb{R} \longrightarrow
    G/K \longrightarrow \Diff_{+}(S^1)
    \longrightarrow 1
  \end{equation}
  is a central $\mathbb{R}$-extension.
\end{proposition}

\begin{pf}
  Note that the equality $h_1 K = h_2 K$ is equivalent to
  $\flux_{\mathbb{R}} (h_1) = \flux_{\mathbb{R}} (h_2)$ for every
  $h_1,h_2 \in G_{\rel}$.
  Since $\flux_{\mathbb{R}}(hgh^{-1}) = \flux_{\mathbb{R}}(g)$ for
  $g \in G$
  and
  $h \in G_{\rel}$,
  we obtain
  \[
    hgK = hgh^{-1}hK = hgh^{-1}K \cdot hK = gK \cdot hK = ghK.
  \]
  The above equality implies that the sequence (\ref{rel_flux_ext})
  is a central extension.
\end{pf}

\begin{definition}
  The central $\mathbb{R}$-extension
  \[
    0 \longrightarrow \mathbb{R} \longrightarrow
    G/K \longrightarrow \Diff_{+}(S^1)
    \longrightarrow 1
  \]
  is called the \textit{flux extension}.
\end{definition}

The $1$-cochain $\tau$
of (\ref{def_eq_tau}) induces a
connection cochain
$\overline{\tau} \in C^1(G/K;\mathbb{R})$
of the flux extension
because of the Corollary \ref{cor_of_flux_and_tau} i).
Thus the formula (\ref{eq_delta_tau})
gives the following Euler cocycle
of flux extension
\[
  -\delta \overline{\tau}(gK,hK)
  =
  \frac{1}{2}
  (\varphi \psi (0) - \varphi(0) - \psi(0))
\]
where $\varphi, \psi$ are defined as in
Proposition \ref{prop_delta_tau} for $g,h \in G$.
Together with the formula (\ref{univ_Euler_cocycle}),
we obtain
\[
  -\delta \overline{\tau}(gK,hK)
  =
  \pi \chi(\mu, \nu)
\]
where $\mu = g|_{\partial D}$ and $\nu = h|_{\partial D}$.
Consequently we have the following:
\begin{theorem}\label{main_thm_1}
  The class $e(G/K)$ is equal to $e_{\mathbb{R}}$
  up to constant multiple.
  More precisely,
  the following holds;
  \[
    e(G/K) = \pi e_{\mathbb{R}}.
  \]
\end{theorem}

\section{Ismagilov-Losik-Michor's cocycle}

Let $M$ be a connected symplectic manifold with an
exact symplectic form $\omega = d\eta$.
Assume that the first Betti number of $M$ is $0$.
Then, there is a $2$-cocycle $C_{\eta,x_0}$
in $C^2(\Symp(M);\mathbb{R})$
introduced by Ismagilov, Losik, and Michor\,\cite{IsLoMi06},
which we call the ILM cocycle.
For $x_0 \in M$,
the ILM cocycle is defined as
\begin{equation}
  C_{\eta,x_0}(g,h)
  = - \int_{x_0}^{h(x_0)} \delta \eta (g)
  = \int_{x_0}^{h(x_0)}\eta^g - \eta
\end{equation}
for $g,h \in \Symp(M)$.
They proves that the cohomology class
$[C_{\eta,x_0}]$ 
is independent of the choice of $x_0 \in M$ and
the potential $\eta \in \Omega^1(M)$.

\begin{remark}
  Let $\mathbb{H}$ be the upper half-plane with a symplectic
  form $dx\wedge dy/y^2$.
  Then the $PSL(2, \mathbb{R})$ is considered as a subgroup of
  $\Symp(\mathbb{H})$
  and
  the ILM cocycle on the group $PSL(2, \mathbb{R})$
  is cohomologous to the area $2$-cocycle
  which is defined as the area of geodesic triangle
  (see \cite{IsLoMi06}).
\end{remark}

The ILM cocycle corresponds to the symplectic form $\omega$ in the
double complex
$(C^{\bullet}(\Symp(M);\Omega^{\bullet}(M)),\delta,d)$
as follows:
Considering a real number as a constant function on $M$,
there is a inclusion
\[
  C^{\bullet}(\Symp(M);\mathbb{R}) \hookrightarrow
  C^{\bullet}(\Symp(M);\Omega^0(M)).
\]
There exists a function $\mathcal{K}_{\eta, x_0}(g)$
satisfying
\[
  d\mathcal{K}_{\eta, x_0}(g) = \delta \eta (g)
  \hspace{2mm}
  \text{ and }
  \hspace{2mm}
  \mathcal{K}_{\eta, x_0}(g)(x_0) = 0
\]
because of $H^1(M;\mathbb{R}) = 0$.
Since the manifold $M$ is connected,
this function $\mathcal{K}_{\eta, x_0}(g)$
is uniquely determined
for each $g \in \Symp(M)$.
Then we obtain the
$1$-cochain
$\mathcal{K}_{\eta, x_0} \in C^1(\Symp(M);\Omega^0(M))$.
Then these cochains $\omega, \eta, \mathcal{K}_{\eta,x_0}$
and $C_{\eta, x_0}$ are connected via the following diagram:
\begin{equation}\label{zigzag}
\vcenter{
\xymatrix{
\omega
\ar@{|->}[r]^-{\delta} &
0 & \\
\eta
\ar@{|->}[u]_-{d} \ar@{|->}[r]^-{\delta}&
\bullet
\ar@{|->}[u]_-{d} \ar@{|->}[r]^-{\delta}&
0\\
&
\mathcal{K}_{\eta,x_0}
\ar@{|->}[u]_-{d} \ar@{|->}[r]^-{\delta}&
\bullet
\ar@{|->}[u]_-{d}\\
& &
-C_{\eta,x_0}
\ar@{|->}[u]_{i}
}
}
\end{equation}

\begin{remark}
  We consider the differential forms
  $\Omega^{\bullet}(M)$ as a right $\Symp(M)$-module
  by pullback.
\end{remark}

Let us return to the disk case.
We consider the ILM cocycle $C_{\eta,x_0}$ on the
unit disk $D$.
Set $H = \Symp(D)$ and
let $i:G \hookrightarrow H$ be the inclusion.
By the equality (\ref{eq_delta_tau_and_ILM}),
the pullback of $C_{\eta,x_0}$ by $i$
is a coboundary, that is,
\begin{equation}\label{tau_and_ILM}
  -\delta \tau = i^* C_{\eta,x_0} \in
  C^2(G;\mathbb{R}).
\end{equation}
So we have proved;
\begin{theorem}\label{main_thm_2}
  On the sequence
  \[
    1 \longrightarrow G_{\rel} \longrightarrow
    G \longrightarrow \Diff_{+}(S^1)
    \longrightarrow 1,
  \]
  the ILM cocycle $C_{\eta,x_0}$ is basic
  and the basic cocycle represents
  $\pi e_{\mathbb{R}}$.
\end{theorem}

\begin{corollary}\label{remofquasimor}
  The pullback $i^* C_{\eta,x_0}$ is a bounded cocycle.
\end{corollary}

\begin{pf}
  Since $\tau$ is a quasi-morphism,
  the boundedness follows.
\end{pf}

\section{The Calabi invariant and extension}

In this section, we investigate the relation
between the ILM cocycle and the Calabi invariant.
First we summarize the results in \cite{Mori16},
which are needed to us.
Let us recall the Calabi invariant and the Calabi extension.
Denote by $H_{\rel}$ the group
consists of the relative symplectomorphism,
that is, the restriction to the boundary $\partial D$
is equal to the identity $\id_{\partial D}$;
\[
  H_{\rel}
  =
  \{ g \in H \mid g|_{\partial D} = \id \}.
\]
Then there is an exact sequence
\[
  1 \longrightarrow H_{\rel} \longrightarrow
  H \longrightarrow \Diff_{+}(S^1)
  \longrightarrow 1.
\]
On the group $H_{\rel}$,
the homomorphism
$\Cal : H_{\rel} \to \mathbb{R}$
is defined by
\begin{equation}\label{def_eq_cal}
  \Cal(g)
  =
  \int_{D} \eta \cup \delta \eta (g)
  =
  \int_{D} \eta^g \wedge (\eta - \eta^g)
\end{equation}
where $\cup$ is the cup product of group cochains.
This homomorphism is called the Calabi invariant
(see \cite{mcduff_salamon}
for details).
It is known that the Calabi invariant
is a surjective homomorphism and
$L = \Ker \Cal$ is
a normal subgroup of $H$.
Dividing the above exact sequence by $L$,
we have the following central extension
\[
  0 \longrightarrow \mathbb{R} \longrightarrow
  H/L \longrightarrow \Diff_{+}(S^1)
  \longrightarrow 1
\]
called the Calabi extension.

\begin{remark}
  There are no inclusion relation between $L$ and
  $K = \Ker \flux_{\mathbb{R}}$ even if we consider them in $G$.
\end{remark}

Next we consider the resulting class $e(H/L)$
of the Calabi extension.
The right-hand side of the formula (\ref{def_eq_cal})
also defines the $1$-cochain
$\tau_0: H \to \mathbb{R}$.
The $1$-cochain $\tau_0 \in C^1(H;\mathbb{R})$
induces a connection cochain
\[
  \overline{\tau_0} \in C^1(H/L;\mathbb{R})
\]
of the Calabi extension.
The curvature $\delta \overline{\tau_0}$
is basic and
this basic cocycle, denoted by the same letter
$\delta \overline{\tau_0}$,
coincides with
$\pi^2 \chi + \pi^2/2$.
Thus we have $e(H/L) = -\pi^2 e_{\mathbb{R}}$
(see Moriyoshi\,\cite[Theorem 2]{Mori16}).

The ILM cocycle also gives a curvature of the
Calabi extension.
Define a $1$-cochain $\kappa \in C^1(H;\mathbb{R})$
by
\[
  \kappa
  =
  \int_{\partial D} \mathcal{K}_{\eta,x_0} \cup \eta.
\]

\begin{remark}
  This $1$-cochain $\kappa$ is a non-zero function.
\end{remark}

\begin{proposition}\label{ilm_is_coboundary}
  Set
  $\tau' = \tau_0 + \kappa
  \in C^1(H;\mathbb{R})$.
  The following holds:
  \[
    -\delta \tau'
    =
    \pi C_{\eta,x_0}.
  \]
\end{proposition}

\begin{pf}
  For $g \in H$, we have
  \begin{align}
  \label{tau_for_cal_by_Kappa}
    \tau_0(g)
    &=
    \int_D \eta^g \wedge \eta \nonumber
    =
    -\int_D \delta \eta (g) \wedge \eta
    =
    -\int_D (d\mathcal{K}_{\eta,x_0}(g))
    \wedge \eta \nonumber \\
    &=
    \int_D \mathcal{K}_{\eta,x_0}(g) \omega
    -
    \int_{\partial D} \mathcal{K}_{\eta,x_0}(g) \eta
    =
    \left( \int_D \mathcal{K}_{\eta,x_0} \cup \omega
    -
    \kappa \right) (g)
  \end{align}
  where the fourth equality follows from the Stokes formula.
  Recalling that $\delta \mathcal{K}_{\eta,x_0} = -C_{\eta,x_0}$,
  we get
  \[
    \delta \tau_0
    =
    \int_D (-C_{\eta,x_0})\cup \omega - \delta \kappa
    =
    -\pi C_{\eta,x_0} - \delta \kappa.
  \]
  This proves $-\delta \tau' = \pi C_{\eta,x_0}$.
\end{pf}

\begin{remark}
  Combining equality (\ref{tau_and_ILM}) and
  Proposition \ref{ilm_is_coboundary},
  we have $\delta (\pi \tau - \tau'|_{G}) = 0$.
  So the function $\pi \tau - \tau'|_{G} : G \to \mathbb{R}$
  is a homomorphism.
  Furthermore, this homomorphism is surjective.
\end{remark}

Next, we show $\tau'$ also gives rise to a connection cochain.

\begin{lemma}\label{properties_kappa}
  For $g \in H, h \in H_{\rel}$,
  we have the following;
  \begin{enumerate}
    \item $\kappa(h) = 0$
    \item $\kappa(gh) = \kappa(g) = \kappa(hg)$.
  \end{enumerate}
\end{lemma}

\begin{pf}
  i) The $1$-cochain $\mathcal{K}_{\eta,x_0}$ satisfies
  \begin{equation}\label{kappa_h}
    d\mathcal{K}_{\eta, x_0}(h) = \delta \eta (h)
    \hspace{2mm}
    \text{ and }
    \hspace{2mm}
    \mathcal{K}_{\eta, x_0}(h)(x_0) = 0
  \end{equation}
  by the definition of $\mathcal{K}_{\eta,x_0}$.
  The restriction of closed form
  $\delta \eta (h)$ to $\partial D$
  is equal to $0$ because the restriction $h|_{\partial D}$
  is the identity.
  There exists a $0$-form $f :D \to \mathbb{R}$
  such that
  \begin{equation}\label{kappa_h_2}
    \delta \eta (h) = df
    \hspace{2mm}
    \text{ and }
    \hspace{2mm}
    f|_{\partial D} = 0.
  \end{equation}
  because of $H^1(D, \partial D; \mathbb{R}) = 0$.
  From the equalities in (\ref{kappa_h}) and (\ref{kappa_h_2}),
  we deduce
  \[
    \mathcal{K}_{\eta,x_0}(h) - f = c
  \]
  where $c$ is a constant in $\mathbb{R}$.
  Evaluating $\mathcal{K}_{\eta, x_0}(h) - f$ at $x_0$,
  we have $c = 0$.
  So we obtain
  \[
    \kappa(h)
    =
    \int_{\partial D} \mathcal{K}_{\eta,x_0} \cup \eta(h)
    =
    \int_{\partial D} f\cdot\eta
    =
    0.
  \]

  ii)
  Note that
  \[
    d(\mathcal{K}_{\eta,x_0}(gh) - \mathcal{K}_{\eta,x_0}(g))
    =
    \eta^g - \eta^{gh}
    =
    d\mathcal{K}_{\eta^g,x_0}(h)
  \]
  and
  $(\mathcal{K}_{\eta,x_0}(gh) -
  \mathcal{K}_{\eta,x_0}(g))(h)(x_0)
  =
  0
  =
  \mathcal{K}_{\eta^g,x_0}(x_0)$.
  By the same argument in i), the restriction
  $\mathcal{K}_{\eta^g,x_0}(x_0)|_{\partial D}$ is equal to $0$.
  So we have
  \[
    \mathcal{K}_{\eta,x_0}(gh)|_{\partial D}
    =
    \mathcal{K}_{\eta,x_0}(g)(x_0)|_{\partial D}
  \]
  and this equality induces $\kappa(gh) = \kappa(g)$.
  Applying the similar argument for
  \[
    d(\mathcal{K}_{\eta,x_0}(hg) - \mathcal{K}_{\eta,x_0}(g))
    =
    \eta^g - \eta^{hg}
    =
    d(\mathcal{K}_{\eta,x_0}(h)^g),
  \]
  we obtain $\kappa(hg) = \kappa(g)$.
\end{pf}

\begin{lemma}\label{tau_ilm_conn}
  For $g \in H, h \in H_{\rel}$,
  the following hold:
  \begin{align}\label{tau_cal_kappa}
    \tau'(gh)
    =
    \tau'(g) + \Cal(h)
    =
    \tau'(hg).
  \end{align}
\end{lemma}

\begin{pf}
  The $1$-cochain $\tau_0$ satisfies
  \[
    \tau_0(gh)
    =
    \tau_0(g) + \Cal(h)
    =
    \tau_0(hg)
  \]
  where $g \in H$ and $h \in H_{\rel}$
  (see \cite[Proposition 4]{Mori16}).
  Combining the above equalities and
  Lemma \ref{properties_kappa}, we obtain
  the equalities (\ref{tau_cal_kappa}).
\end{pf}

From Lemma \ref{tau_ilm_conn},
the cochain $\tau'$ induces the connection cochain
$\overline{\tau'} \in C^1(H/L;\mathbb{R})$
of Calabi extension.
So $-\delta \overline{\tau'}$
gives the Euler class
$e(H/L) \in H^2(\Diff_{+}(S^1);\mathbb{R})$.
Recalling $e(H/L) = -\pi^2 e_{\mathbb{R}}$,
we obtain:
\begin{theorem}\label{main_thm_3}
  The ILM cocycle $C_{\eta,x_0}$ is basic
  with respect to the exact sequence
  \[
    1 \longrightarrow H_{\rel} \longrightarrow
    H \longrightarrow \Diff_{+}(S^1)
    \longrightarrow 1
  \]
  and the cohomology class of
  basic cocycle
  is equal to $\pi e_{\mathbb{R}}$.
\end{theorem}

\begin{remark}
  The integration in definition of ILM cocycle
  depends only on the restrictions of $g, h$
  to $\partial D$.
  Hence, from Corollary \ref{remofquasimor},
  the ILM cocycle $C_{\eta,x_0}$ is also a bounded cocycle.
  Since $\delta \tau'$
  is equal to $-\pi C_{\eta,x_0}$,
  the cochain $\tau'$ gives rise to a
  quasi-morphism on $H$.
  In \cite{Mori16},
  it is proved that the cochain
  $\tau_0$ is a quasi-morphism.
  Thus the cochain $\kappa$ is also a quasi-morphism
  on $H$
  because of $\tau' = \tau_0 + \kappa$.
\end{remark}

\end{document}